\documentclass[10pt]{amsart}
\usepackage{amsmath,amssymb,amscd,enumerate,bbm}
\usepackage{mathrsfs}
\usepackage{cite}

\renewcommand{\MR}[1]{}

\newcommand{\titl}{GAMMA STRUCTURES AND GAUSS'S CONTIGUITY}

\title{{\titl}}
\author{V. Golyshev and A. Mellit}
\date{}

\newcommand{\cal}{\mathcal}

\def\I{{i}}
\def\P{{\Bbb{P}}}
\def\R{{\Bbb{R}}}
\def\C{{\Bbb{C}}}
\def\Q{{\Bbb{Q}}}
\def\Z{{\Bbb{Z}}}

\def\HH{{\cal{H}}}

\def \DD {\mathcal{D}}

\def\OO{{\mathcal{O}}}
\renewcommand{\phi}{{\varphi}}
\def\FF{{\cal{F}}}
\def\GG{{\cal{G}}}
\def\Alpha{{\mathrm{A}}}
\def\Beta{{\mathrm{B}}}
\newcommand{\gammastructure}{\mathbf{[\Gamma]}}

\newcommand\Exp[1]{e^{#1}}

\newcommand{\transpose}{{\mathrm{t}}}
\newcommand{\Rep}{{\mathrm{Rep}}}
\newcommand{\Id}{\mathrm{Id}}

\newcommand{\cD}{{\cal D}}

\newcommand{\FT}{{\mathrm{FT}}}
\newcommand{\VM}{{\mathrm{V}}}

\newcommand{\Spec}{{\text{Spec }}}

\newcommand{\Real}{\mathop{\mathrm{Re}}}
\newcommand{\Imaginary}{\mathop{\mathrm{Im}}}

\newcommand{\Arg}{\mathop{\mathrm{arg}}}

\newcounter{pphcounter}[section]
\renewcommand{\thepphcounter}{\thesection.\arabic{pphcounter}}
\newcommand{\pph}{\bigskip \refstepcounter{pphcounter}
    \bf  \thepphcounter. \rm}

\makeatletter
\@addtoreset{equation}{pphcounter}
\makeatother

\newcommand{\coro}{\bf Corollary. \rm}
\newcommand{\defi}{\bf Definition. \rm}
\newcommand{\rema}{\bf Remark. \rm}
\newcommand{\propo}{\bf Proposition. \rm}
\newcommand{\theo}{\bf Theorem. \rm}

\renewcommand{\proof}{\bigskip \bf Proof. \rm}

\def\from{:}

\def\A1{{{\Bbb{A}}^1}}
\def\P1{{{\Bbb{P}}^1}}

\def\Gm{{\bf G_m}}

\def\SL2{{\mathrm SL2}}

\renewcommand{\P}{{\Bbb{P}}}

\newcommand{\benum}{\begin{enumerate}}
\newcommand{\eenum}{\end{enumerate}}

\newcommand{\sols}{f}
\newcommand{\mystrut}{{\vphantom{\vrule height 1.2em}}}

\begin{document}

\begin{center}
\maketitle

\bigskip

\bigskip

\end{center}

\bigskip

\parbox{380pt}{\small \bf Abstract. \rm
We introduce gamma structures on regular hypergeometric D--modules
in dimension $1$
as special one--parametric systems of solutions on the compact subtorus.
We note that a balanced gamma product
is in the Paley--Wiener class and show that the monodromy
with respect to the gamma structure is expressed algebraically
in terms of the hypergeometric exponents. We
compute the hypergeometric monodromy explicitly
in terms of certain diagonal matrices,
Vandermonde matrices and their inverses (or generalizations of those in the resonant
case).
}

\bigskip
\bigskip
\bigskip

A hypergeometric D--module with rational indices is motivic, i. e.
may be realized as a constituent of the pushforward of the constant
D--module $\OO$ in a pencil of varieties over $\Gm$ defined over
$\overline{\mathstrut \Q}$. The de Rham to Betti comparison arises in
each fiber; as a result, the vector space of solutions is endowed
with two $K$--rational structures for a number  field $K$. On the other
hand, no rational structure exists in the case of irrational exponents,
and yet one still wishes to have the benefits of the Dwork/Boyarsky method
of parametric exponents.

A substitute is the \emph{gamma structure} on a hypergeometric D--module
which
manifests itself as a rational
structure in the case of rational exponents and gives rise to an extension
of the Betti to de Rham comparison in the non-motivic direction.

Operating in this framework,
one might try to study period matrices of
traditional motives by representing them as limiting cases of hypergeometric ones,
or even
degenerate the hypergeometric period matrix into a resonant singularity.
Reverting this process yields a perturbation of a Tate type period to an expression  in gamma--values, cf \cite{Golyshev08a}.

F. Baldassarri has emphasized that the key to hypergeometric monodromy
 is Gauss's contiguity principle:
 with a translation of the set of indices
 \footnote{We adopt the terminology in which, in
 $x^\alpha$, $\alpha$ is the index, and $\exp( 2 \pi i \alpha)$, the exponent.}
 by a vector in an integral lattice is
 associated an explicit isomorphism of the respective D--modules,
 whose shape leads one to an a priori guess on the shape of the monodromy.
 Y. Andre has remarked that the situation is even better with $p$--adic hypergeometrics,
 as the translation lattice is dense in the space of indices.
 We introduce the gamma structure and  replication
  as a means to make up for the lack of density
 of the translations in the complex case by interpolating the shifts to
 non--integral ones.

\bigskip
\begin{center}
------
\end{center}
\bigskip

\bigskip

We follow Katz's treatment \cite{Katz90} of hypergeometrics in order to fix our basics .
Let $\Gm=\Spec\C[z,z^{-1}]$ be a one--dimensional torus. By
$\cD$ denote the algebra of differential operators on $\Gm$, by
$D$~denote the differential operator $z\dfrac{\partial}{\partial
z}$. One has $\cD=\C[z,z^{-1},D]$.

 \defi
Let $n$ and
$m$~be a pair of nonnegative integers.
Let  $P$ and $Q$ be polynomials
of degrees $n$ and $m$ respectively. Define {\it a hypergeometric
differential operator\/} of type $(n,m)$:
$$
\HH (P,Q)=P(D)-zQ(D),
$$
and {\it the hypergeometric  $D$-module\/}:
$$
\HH (P,Q) = \mathscr{D}/ \mathscr{D}\HH(P,Q).
$$

If $P(t)=p \prod_i(t-a_i)$, \ $Q(t)=q \prod_j(t-b_j)$, \
$\lambda=p/q$, we shall write
$$
\HH_\lambda(a_i,b_j) = \DD /\DD \left(\lambda \prod_i(D-a_i)-
z\prod_j(D-b_j)\right)).
$$

 \propo
\benum
\item
$\HH_\lambda(a_i,b_j)$ is an irreducible
$D$-module on $\Gm$ if and only if  $P$ and $Q$  have  no common
zeros $\mod \Z$.


\item
Let
$\HH_\lambda(a_i,b_j)$~be an irreducible hypergeometric
$D$-module. If $\#\{a_i\} \ne \#\{b_j\}$ then
$\HH_\lambda(a_i,b_j)$ is a differential equation on $\Gm$. If
$\# \{a_i\}=\# \{b_j\}$, put $U=\Gm\setminus\{\lambda\}$ and
let  $j\from U \to \Gm$ be the respective open immersion. Then
$\HH_\lambda(a_i,b_j)$  is a differential equation on $U$, \
$\HH_\lambda(a_i,b_j)=j_{!*}j^*\HH_\lambda(a_i,b_j)$, and the
local monodromy of its solutions around $\lambda$ is a
pseudoreflection.

\item
Fix a $\lambda$.
Let $\HH_\lambda(a_i,b_j)=\HH(P,Q)$ be an irreducible $D$-module
on $\Gm$. Then the isomorphism class of \linebreak
$\HH_\lambda(a_i,b_j)$ depends only on the sets $\{ a_i \mod \Z
\} $ and $\{ b_j \mod \Z \}$.

\item
Let
$\HH_\lambda(a_i,b_j)=\HH(P,Q)$ be an irreducible $D$-module on
$\Gm$ of type $(n,m)$. If $n \ge m$ respectively, $m \ge
n$), then the eigenvalues of the local monodromy at the
regular singularity $0$ are $\exp(2\pi i a)_{P(a)=0}$ (resp.,
the eigenvalues of the local monodromy at the regular singularity
$\infty$ are $\exp(2\pi i b)_{Q(b)=0}$); to each eigenvalue
of the local monodromy at $0$ (resp. at $\infty$)
corresponds the unique Jordan block.


\item
The isomorphism
class $\HH_\lambda(a_i,b_j)$ determines the type $(n,m)$, the sets
$\{a_i \mod \Z\}$, $\{b_j \mod \Z\}$ with multiplicities and,
in the $n=m$ case, the scalar $\lambda$.

\item
Let $\FF$, $\GG$ be two irreducible local
systems on  $(\Gm\setminus\{\lambda\})^{an}$ of the same rank
$n \ge 1$. Assume that:

\benum

\item the local monodromies of both systems at $\lambda$ are
pseudoreflections;

\item the characteristic polynomials of the local systems
$\FF$ and $\GG$ at $0$ are equal;

\item  the characteristic polynomials of the local systems
$\FF$ and $\GG$ at $\infty$ are equal.

\eenum

Then there is an isomorphism  $\FF\cong \GG$.

\eenum
\bigskip

%
%


%
%
%
%

\begin{center}
------
\end{center}
\bigskip

We say that a holomorphic function $ {\Phi} (s)$
is of Paley--Wiener type if it is a
Fourier transform of a function/distribution $H$ on $\R$
with compact support.
Assume that $f$ satisfies a linear homogeneous recurrence $\mathrm{R}$
with polynomial coefficients.
Then, for any periodic distribution $p(s)$, the product $p(s){\Phi}(s) $
satisfies $\mathrm{R}$ as well, its inverse Fourier transform being a solution
to the DE that is the inverse Fourier transform of $\mathrm{R}$.
In particular, let $p= \Delta^{t}= \sum_{l \in \Z} \delta_{t+l}.$
We thus get a system of solutions $S_{t}$.

Define now a \emph{balanced gamma product} by
$$ \mathbf{\Gamma} (s)
= \dfrac{1} {\mystrut \prod_{i=1}^n \Gamma(s-\alpha_i+1)
\prod_{j=1}^n \Gamma(-s + \beta_j + 1)}.$$ Then
$ \mathbf{\Gamma} (s)$ is of Paley--Wiener type, and the considerations above apply. The
corresponding hypergeometric equation is $\HH_{\lambda}(\alpha_i, \beta_i)$ with $\lambda = (-1)^n$.
The inverse Fourier transform $h$ of $\mathbf{\Gamma} (s)$ is a solution of the hypergeometric
equation on the universal cover of the unit circle. This solution is supported on
$[-\dfrac{n}2, \dfrac{n}2]$, which is a union of $n$ segments of length $1$, each segment
being identified with the unit circle without the singular point. Thus we obtain $n$ solutions
of the equation on the unit circle, which form a basis, which we denote by $\sols$,
and a one--parameter family of solutions $S_{t}$.

Looked at from this viewpoint, the non--resonant hypergeometric monodromy
(i.e. one with distinct $\alpha$'s
and $\beta$'s $\mod \Z$) can be computed easily as follows.
Construct a basis of solutions given by the power series in the neighborhood
of $0$, and similarly in the neighborhood of $\infty$.
In the notation adopted above, the former basis is simply $\{ S_{\alpha_k} \}$, and the latter,
$\{ S_{\beta_{k'}} \}$. It is clear that
$$S_{\alpha_k} = \sum {\sols_m} \exp (2 \pi i ([n/2]-m) \alpha_k) \text{ and }
S_{\beta_{k'}} = \sum {\sols_m} \exp (2 \pi i ([n/2]-m) \beta_{k'}).$$
In the basis $S_\alpha$ the monodromy around $0$ is diagonal with eigenvalues $\exp (2 \pi i \alpha_k)$,
in the basis $S_\beta$ the monodromy around $\infty$ is diagonal with eigenvalues $\exp (2 \pi i \beta_{k'})$,
and the relation above of each basis to $\sols$ is a means to glue up the two. The present paper
is an elaboration of this concise argument in a possibly resonant case.

Put $\partial=\frac{1}{2\pi\I}\frac{\partial}{ \partial t}$. We take the
liberty of denoting this derivation also as $\frac{\partial}{2\pi\I \partial t}$. In a resonant case, one considers derived periodic distributions $p=(-1)^r \partial^r\Delta^{t}$.
A \emph{replication} of $h$ is the resulting inverse Fourier transform of $\mathbf\Gamma \cdot p$. The \emph{gamma structure} $\gammastructure$ on $\HH_\lambda(\alpha, \beta)$ associated with $\mathbf\Gamma$
is defined to be the set of all replications of $h$.
Theorem \ref{theorem_on five--tuples} states that the hypergeometric monodromy expressed in terms of the basis of local solutions at $0$ (resp. $\infty$) that are in the gamma structure is given by products of generalized diagonal matrices and Vandermonde
matrices (and their inverses), whose entries depend algebraically on the hypergeometric exponents, cf \cite{Levelt61}.

\bigskip
\bigskip
\section{The Paley-Wiener property of gamma products}
\pph \label{estimate} \propo There exists $C>0$ such that for all $s\in\C$
$$
\bigg\lvert{\dfrac{1}{\Gamma(s)}}\bigg\rvert < C (1 + |s|)^{\frac12 - \Real s}\; \Exp{\Arg{s} \Imaginary s + \Real s},
$$
where $\Arg{s}$ is chosen to be in $[-\pi, \pi]$.

\proof Apply Stirling's approximation \cite[13.6]{WW27}
$$
\Gamma(s)=\sqrt{\frac{2\pi}s} \left(\frac{s}{e}\right)^s (1 + o(1)),
$$
which holds when $\Real s \ge 0,\, \mid s \mid \to \infty$. The proof of the estimate in this case is straightforward. When $\Real s \le 0$ we apply
$$
\frac{1}{\Gamma(s)} = \frac{(-s)\Gamma(-s)(\Exp{\pi\I s} - \Exp{-\pi\I s})}{2\pi\I}.
$$
We may assume $\Imaginary s \ge 0$ without loss of generality. Applying Stirling's approximation and noting that $\Arg{(-s)} = \Arg s - \pi$ gives the proof in this case.

\pph \defi Let $PW_R$ be the space of entire functions $f$ such that for some $C, \mu\in\R$ the following estimate holds:
$$
|f(s)| < C (1 + |s|)^\mu \Exp{2 \pi R |\Imaginary s|} \quad \text{for all $s\in\C$}.
$$

\pph {\bf Fourier transform.} The Fourier transform is a continuous automorphism of the space of tempered distributions $S'(\R)$. We recall that a tempered distribution is a continuous linear functional on the space of infinitely differentiable functions of rapid decay. The Fourier transform is given on integrable functions by
$$
(\FT f)(s) = \int_{\R} \Exp{-2\pi\I\phi s} f(\phi) d\phi.
$$

\pph {\bf The Paley-Wiener theorem.} \cite[7.23]{Rudin73} A tempered distribution is the Fourier transform of a distribution supported on $[-R, R]$ if and only if it can be extended to a function in $PW_R$.

\bigskip
We could not find a reference for the following classical--looking theorem.

\pph \theo For $\alpha_1, \ldots, \alpha_n, \beta_1, \ldots, \beta_n \in \C$ 
the following function is in the space $PW_{\frac{n}2}$:
$$
\mathbf{\Gamma}_{\alpha, \beta}(s) \,:=\, \dfrac{1}{\mystrut \prod_{i=1}^n \Gamma(s - \alpha_i + 1) \prod_{i=1}^n \Gamma (-s + \beta_i + 1)} \quad (s\in\C).
$$

\proof Apply Proposition \ref{estimate}. \qed

\medskip
We will construct the inverse Fourier transform of $\mathbf \Gamma_{\alpha, \beta}(s)$ explicitly in \ref{pph:h_explicit}.

\bigskip
\bigskip
\section{Vandermonde matrices, diagonal matrices, cyclic matrices}
Let $\Alpha = (\Alpha_1,\ldots,\Alpha_{n_A})$ be a tuple of distinct non-zero complex numbers. Let $m_\Alpha = (m_{\Alpha_1},\ldots,m_{\Alpha_n})$ with $m_{\Alpha_j}\geq 0$ be integers, which we call {\em multiplicities}. We will consider square matrices of size $n=\sum_j m_{\Alpha_j}$ of three types.

\pph {\bf Generalized Vandermonde matrix.} For any $l\in\Z$ this is the $n\times n$ matrix, denoted $\VM_{\Alpha, m_\Alpha, l}$, whose rows are indexed by pairs $(j, r)$ with $j=1,\ldots,n_\Alpha$ and $r=0,\ldots,m_{\Alpha_j}-1$, columns are indexed by $k=0,\ldots,n-1$, and elements are
$$
(\VM_{\Alpha, m_\Alpha, l})_{(j,r),\,k} = (l-k)^r A_j^{l-k}\quad \text{where $0\leq k < n$, $1\leq j \leq n_\Alpha$, $0\leq r < m_{\Alpha_j}$, so that}
$$
$$
\VM_{\Alpha, m_\Alpha, l} = \begin{pmatrix}
A_1^l & A_1^{l-1} & \cdots & A_1^{l-n+1} \\
l A_1^l & (l-1) A_1^{l-1} & \cdots & (l-k+1) A_1^{l-n+1} \\
\vdots & \vdots & & \vdots
\end{pmatrix}.
$$

\pph {\bf A block-diagonal matrix.} Let $D_{\Alpha, m_\Alpha}$ be the $n\times n$ matrix whose rows and columns are indexed by pairs $(j, r)$ as above, and elements are
$$
(D_{\Alpha, m_\Alpha})_{(j,r),\,(j', r')} =
\begin{cases}
\binom{r}{r'} \Alpha_j & \text{if}\; j=j',\; r'\leq r, \\
0 & \text{otherwise}.
\end{cases}
$$
The matrix $D_{\Alpha, m_\Alpha}$ is lower-triangular with diagonal elements $A_j$.

\pph \propo The matrix $\VM_{\Alpha, m_\Alpha, l}^{-1} D_{\Alpha, m_\Alpha} \VM_{\Alpha, m_\Alpha, l}$ is of the {\em cyclic form}, i.e.
$$
\VM_{\Alpha, m_\Alpha, l}^{-1} D_{\Alpha, m_\Alpha} \VM_{\Alpha, m_\Alpha, l} =
\begin{pmatrix}
* & 1 & 0 & \cdots & 0\\
* & 0 & 1 & \cdots & 0\\
\vdots & \vdots & \vdots & & \vdots \\
* & 0 & 0 & \cdots & 1 \\
* & 0 & 0 & \cdots & 0
\end{pmatrix}.
$$

\proof
Let $P(x) = \prod_{j=1}^{n_\Alpha} (x - \Alpha_j)^{m_{\Alpha_j}}$ and
consider the finite algebra $R=\C[x, x^{-1}]/P(x)$. Put $\bar x = x \mod P(x)$,
 and choose $\bar x^l, \bar x^{l-1},\ldots,\bar x^{l-n+1}$ for a $\C$-basis of $R$.
 Identify the standard $n$--dimensional space $\C ^n$ with the standard basis with the space of principal
 parts of Laurent polynomials in $x$ at $A_1, \dots, A_n$. The matrix  $\VM_{\Alpha, m_\Alpha, l}$ is the matrix of the linear operator $R \to \C ^n$ that maps an element $\bar f \in R$ to a vector with components $((x \frac{d}{dx})^r f)(\Alpha_j)$.
 Multiplication by $x$ is a linear operator on the space of principal parts. By the
 Leibnitz rule, it
 transforms the vector $0,\dots, 1,\dots, 0$ with the single $1$ at the place that corresponds to $j, r'$ into
 $0, \dots, A_j, (r'+1) A_j,\dots, {m_{A_j}-1 \choose r'}A_j, \dots, 0$, or, in other words,
 is given by the matrix $D_{\Alpha, m_\Alpha}$.
The matrix
$\VM_{\Alpha, m_\Alpha, l}^{-1} D_{\Alpha, m_\Alpha} \VM_{\Alpha, m_\Alpha,l}$ is then the matrix of multiplication by $x$ expressed in
the basis $\bar x^l, \bar x^{l-1},\ldots,\bar x^{l-n+1}$, and is therefore
cyclic.

\section{Local solutions}\label{local_sol}
Let $\alpha=(\alpha_1, \ldots, \alpha_n)$, $\beta=(\beta_1, \ldots, \beta_n)$ be tuples of complex
numbers such that $\alpha_i \neq \beta_{i'} \mod \Z$ for all $i, i'$. Put
$$
\mathbf{\Gamma}(s)\,=\, \mathbf{\Gamma}_{\alpha, \beta}(s) \,=\, \dfrac{1}{\mystrut \prod_{i=1}^n \Gamma(s  - \alpha_i + 1) \prod_{i=1}^n \Gamma (-s + \beta_i + 1)} \quad (s\in\C).
$$

\pph {\bf Basis at $0$.} Let $\Alpha_1,\ldots,\Alpha_{n_\Alpha}$ be all distinct values in $\Exp{2\pi\I\alpha_1},\ldots,\Exp{2\pi\I\alpha_n}$, and put $m_{\Alpha_j} = \#\{k:\Exp{2\pi\I\alpha_k} = \Alpha_j\}$. Define a map $\nu:\{\Alpha_1,\ldots,\Alpha_{n_\Alpha}\} \rightarrow \{\alpha_1, \ldots, \alpha_n\}$ by the condition that $\Exp{2\pi\I\nu(\Alpha_j)} = \Alpha_j$ and $\nu(\Alpha_j)$ has the minimal real part
among such $\alpha_i$. Put
$$
S_{\Alpha_j, r}(z) = \sum_{l=0}^\infty \left(\dfrac{\partial}{2\pi\I \partial t}\right)^r \bigg\rvert_{t=\nu(\Alpha_j)} \left( \mathbf\Gamma(l + t)\, z^{l + t} \right)\quad (r<m_{\Alpha_j}).
$$

\pph \propo The formal (log) power series $S_{\Alpha_j, r}$ for $j = 1, \ldots, n_\Alpha$, $r=1,\ldots,m_{\Alpha_j}$ form a basis of solutions at $0$ for the differential equation $\HH_{\lambda}(\alpha_i, \beta_i)$ with $\lambda = (-1)^n$.

\proof First we prove that $S_{\Alpha_j, r}(z)$ is a solution of the differential equation. Note that $\mathbf\Gamma(l + t)$ has $0$ of order $m_j$ at $\alpha_{i(j)}$ if $l<0$, $l\in\Z$. Therefore the sum defining $S_{\Alpha_j, r}(z)$ can be formally replaced with the sum over all $l\in\Z$. The derivation $z\frac{d}{dz}$ acts on the expansion coefficients of these series by sending $f(s)$ to $s f(s)$  Thus the statement follows from the identity for $\mathbf\Gamma$:
$$
\mathbf\Gamma(s) \;\prod_{i=1}^n (s - \alpha_i) \;=\; \mathbf\Gamma(s-1) \;\prod_{i=1}^n (-(s-1) + \beta_i).
$$

To see that $S_{\Alpha_j, r}(z)$ are linearly independent we first note that $S_{\Alpha_j, 0}(z)$ are non-zero because of the condition $\alpha_i \neq \beta_{i'} \mod\Z$. Let $M_0^\Alpha$ be the local monodromy operator at $0$. As will be shown in the next paragraph,
\[
(\Alpha_j^{-1} M_0^\Alpha - \Id)^{r+1} S_{\Alpha_j, r} = 0,\quad (\Alpha_j^{-1} M_0^\Alpha - \Id)^r S_{\Alpha_j, r} = r! S_{\Alpha_j, 0},
\]
and the statement follows.

\pph \propo\label{prop:monodromy_a} The monodromy around $0$ of the basis $S_{\Alpha_j, r}$ is given by the matrix $M_0^\Alpha = D_{\Alpha, m_\Alpha}^\transpose$.

\proof This follows from
\begin{align*}
M_0^\Alpha S_{\Alpha_j, r}(z) &= \sum_{l=0}^\infty \left(\dfrac{\partial}{2\pi\I \partial t}\right)^r \bigg\rvert_{t=\alpha_{i(j)}} \left(\mathbf\Gamma(l + t)\, z^{l + t}  \Exp{2\pi\I t} \right)\\
&= \Alpha_j \sum_{p = 0}^r \binom{r}{p} S_{\Alpha_j, p}.
\end{align*}

\pph {\bf Basis at $\infty$.} Analogously,
let $\Beta_1,\ldots,\Beta_{n_\Beta}$ be all distinct values in $\Exp{2\pi\I\beta_1},\ldots,\Exp{2\pi\I\beta_n}$, and put $m_{\Beta_j} = \#\{k:\Exp{2\pi\I\beta_k} = \Beta_j\}$. Define a map $\nu':\{\Beta_1,\ldots,\Beta_{n_\Beta}\} \rightarrow \{\beta_1, \ldots, \beta_n\}$ by the condition that $\Exp{2\pi\I\nu'(\Beta_j)} = \Beta_j$ and $\nu'(\Beta_j)$ has the maxinal real part
among such $\beta_i$.
We put
$$
S_{\Beta_j, r}(z) = \sum_{l=0}^\infty \left(\dfrac{\partial}{2\pi\I\partial t}\right)^r \bigg\rvert_{t=\beta_{i(j)}} \left( \mathbf\Gamma(-l + t)\, z^{-l + t} \right)\quad (r<m_{\Beta_j}).
$$

\pph \propo The formal power series $S_{\Beta_j, r}$ for $j = 1, \ldots, n_\Beta$, $r=1,\ldots,m_{\Beta_j}$ form a basis of solutions at $\infty$ for the differential equation $\HH_{\lambda}(\alpha_i, \beta_i)$ with $\lambda = (-1)^n$. \qed

\pph \propo \label{prop:monodromy_b} The monodromy around $\infty$ of the basis $S_{\Beta_j, r}$ is given by the matrix $M^\Beta_{\infty} = \left(D_{\Beta, m_\Beta}^\transpose\right)^{-1}$. \qed

\section{Solutions on the unit circle} \label{sec:solutions}
Let $\alpha_i$, $\beta_i$, $\mathbf\Gamma$ be as in the previous section. We construct a basis of solutions of the hypergeometric equation on the unit circle using the Paley-Wiener property of the gamma product. Let $h$ be the distribution supported on $[-\frac{n}2, \frac{n}2]$ whose Fourier transform is $\mathbf\Gamma$. Put
$$
f_k = h|_{(-\frac{n}2 + k, -\frac{n}2 + k + 1)} \quad (k=0,\ldots n-1),
$$
where the interval $(-\dfrac{n}2 + k, -\dfrac{n}2 + k + 1)$ is identified with $S^1 \setminus \{\lambda\}$ in the natural way. As always, $\lambda = (-1)^n$.

\pph \propo The distributions $f_k$ are smooth functions which satisfy the differential equation $\HH_{\lambda}(\alpha_i, \beta_i)$.

\proof
It follows from the formal properties of the Fourier transform that $f_k$ satisfies the differential equation as a distribution. To verify that $f_k$ is a smooth function we will construct $h$ in a different way.

\pph {\bf The case $n=1$.} Suppose for a moment that $n=1$, $\alpha_1 = \alpha$, $\beta_1 = \beta$, $\Real(\beta-\alpha)>0$.
Put
$$
h_{\alpha, \beta} (\phi) \,=\,
\begin{cases}
\dfrac{\Exp{2\pi\I\alpha\phi} (1 + \Exp{2\pi\I\phi})^{\beta-\alpha}}{\Gamma(\beta - \alpha + 1)} & \phi\in\left(-\dfrac12,\dfrac12\right)\\
0 & \phi\notin\left(-\dfrac12,\dfrac12\right)
\end{cases}.
$$

\pph \propo One has $\FT h_{\alpha, \beta} \;=\; \mathbf\Gamma_{\alpha, \beta}$.

\proof Consider the unit disk and cut away the segment $[-1,0]$. Let $S$ be the path on the unit circle which starts and ends at $-1$ and goes  counterclockwise. The integral defining the left hand side can be written as
$$
\dfrac{1}{2\pi\I\, \Gamma(\beta-\alpha+1)} \int_S z^{-s + \alpha} (1 + z)^{\beta-\alpha} \dfrac{d z}{z}.
$$
We deform $S$ to the path which first goes from $-1$ to $0$ just below the cut and then goes from $0$ to $-1$ just above the cut. Making change of variables $z=-u$ the integral becomes
\begin{multline*}
\dfrac{1}{2\pi\I\, \Gamma(\beta-\alpha+1)} (\Exp{\pi\I (-s + \alpha)} - \Exp{\pi\I (s - \alpha)}) \int_0^1 u^{-s + \alpha - 1} (1-u)^{\beta - \alpha} du \\
= \dfrac {\sin{\pi(-s + \alpha)}}{\pi\, \Gamma(\beta-\alpha+1)}\; \Beta(-s + \alpha, \beta - \alpha + 1) = \mathbf\Gamma_{\alpha, \beta} (s).
\end{multline*}

\pph \label{pph:h_explicit} This immediately implies the statement in the case when $\Real(\beta_i-\alpha_i)>0$,
because
$$
h = h_{\alpha_1, \beta_1} * \cdots * h_{\alpha_n, \beta_n},
$$
and each function in the convolution is smooth on $(-\dfrac12,\dfrac12)$, so the only non-smooth points of $h$ are the points $-\dfrac{n}2 + k$, $k=0,\ldots,n$.

\pph To prove the general case we note that $\mathbf\Gamma(s)=\mathbf\Gamma_{\alpha, \beta}(s)$ is the product of a polynomial $R(s)$ and $\mathbf\Gamma'(s) = \mathbf\Gamma_{\alpha, \beta' + m}(s)$ for a positive integer $m$ such that $\Real(\beta_i + m - \alpha_i)>0$ holds. Denoting by $f'_k$ the corresponding solution for $(\alpha_i), (\beta_i + m)$, which is smooth, $f_k(\phi) = R(\dfrac{1}{2\pi\I} \dfrac{d}{d\phi}) f_k'(\phi)$, therefore $f_k$ is smooth as well.

\section{Hypergeometric monodromy}
Let $\alpha_i$, $\beta_i$, $\mathbf\Gamma$, $\Alpha_j$, $\Beta_j$, $h$, $f_k$ be as in the previous sections. Assume $\alpha_i$, $\beta_i$ are real.

\pph {\bf Analytic continuation.}\label{pph:analytic_continuation} The formal series $S_{\Alpha_j, r}$ converge  on the universal cover of the punctured open unit disk and define analytic functions thereon. We will use polar coordinates, $z = \rho \Exp{2\pi\I\phi}$. The restriction of $S_{\Alpha_j, r}$ to the universal cover of the circle of radius $\rho$ is hence given as
$$
S_{\Alpha_j, r}(\rho \Exp{2\pi\I\phi}) \;=\; \sum_{l=0}^\infty \left(\dfrac{\partial}{2\pi\I \partial t}\right)^r \bigg\rvert_{t=\alpha_{i(j)}} \left( \mathbf\Gamma(l + t)\, \rho^{l + t} \Exp{2\pi\I \phi (l + t)} \right).
$$
Since the infinite sum also converges in the topology of $S'$ for a fixed value of $\rho$, the identity above holds in the space $S'$. Therefore we may apply the Fourier transform (from $S'$ to $S'$) termwise. If we denote $S_{\Alpha_j, r, \rho}(\phi) = S_{\Alpha_j, r}(\rho \Exp{2\pi\I\phi})$, we obtain
$$
\FT S_{\Alpha_j, r, \rho} \;=\; \sum_{l=0}^\infty \sum_{p=0}^r (-1)^{p} \binom{r}{p}  \left(\dfrac{\partial}{2\pi\I \partial t}\right)^{r-p} \bigg\rvert_{t=\alpha_{i(j)}} \left( \mathbf\Gamma(l + t)\, \rho^{l + t} \right) \; \delta_{l + \alpha_{i(j)}}^{[p]},
$$
where $\delta_x^{[p]}$ is the $p$~-th derivative of the $\delta$~-distribution concentrated at $x$ divided by $(2\pi\I)^p$. When $\rho$ tends to $1$ from below we obtain the following identity in $S'$:
\begin{align*}
\lim_{\rho\rightarrow 1-} \FT S_{\Alpha_j, r, \rho} \;&=\; \sum_{l=0}^\infty \sum_{p=0}^r (-1)^{p} \binom{r}{p}  \left(\dfrac{\partial}{2\pi\I \partial t}\right)^{r-p} \bigg\rvert_{t=\alpha_{i(j)}} \mathbf\Gamma(l + t) \; \delta_{l + \alpha_{i(j)}}^{[p]} \\
&=\; \mathbf\Gamma \cdot (-1)^r \sum_{l\in\Z} \delta_{l + \alpha_{i(j)}}^{[p]}.
\end{align*}
Note that the product of distributions above is well-defined since $\mathbf\Gamma$ is of Paley-Wiener type.

\pph {\bf Replication of distributions with compact support.}
Let $q \in \C$, $|q| = 1$, $q= \Exp{2\pi\I t}$. Put
$$
\Delta_{q, r} = \sum_{l\in\Z} l^r q^l\delta_l,\quad
\Delta^{q, r} = (-1)^r \sum_{l\in\Z} \delta_{l + t}^{[r]}.
$$
Since $\FT \sum_{l\in\Z} \delta_l = \sum_{l\in\Z} \delta_l$ (the Poisson summation formula), we have $\FT\Delta_{q, r} = \Delta^{q, r}$.

Let $g$ be a distribution on $\R$ with compact support. Its {\em replication of order $r$ and parameter $q$} is a tempered distribution on $\R$ defined as
$$
\Rep_{q, r} g \;:=\; g * \Delta_{q, r}.
$$
One has $\FT\, \Rep_{q, r} g = (\FT g) \cdot \Delta^{q, r}$.

\pph \defi We define the gamma structure $\gammastructure$ on $\HH_{(-1)^n}(\alpha_i, \beta_i)$ associated to $\mathbf\Gamma$ to be the set of all replications of $h$.

\pph \theo \label{thm:replication} Denote by $h$ the inverse Fourier transform of $\mathbf\Gamma$, as in the beginning of Section \ref{sec:solutions}. Then the solutions $S_{\Alpha_j, r}$ around $0$ and $S_{\Beta_j, r}$ around $\infty$ as functions on the universal cover of the unit circle belong to the gamma structure $\gammastructure$. Concretely,
\benum
\item the analytic continuation of the local solution $S_{\Alpha_j, r}$ to the universal cover of the unit circle can be represented as the replication of $h$:
$$
\lim_{\rho\rightarrow 1-} S_{\Alpha_j, r, \rho} = \Rep_{\Alpha_j, r} h;
$$
\item the analytic continuation of the local solution $S_{\Beta_j, r}$ to the universal cover of the unit circle can be represented as the replication of $h$:
$$
\lim_{\rho\rightarrow 1+} S_{\Beta_j, r, \rho} = \Rep_{\Beta_j, r} h;
$$
\eenum

\proof The inverse Fourier transform is continuous on $S'$ \cite[7.15]{Rudin73}. Hence the theorem follows from the limit formula in \ref{pph:analytic_continuation}.

\pph Restricting both sides of the formula in Theorem \ref{thm:replication} ({\rm i}) to the interval $(-\dfrac{n}2+l,-\dfrac{n}2+l+1)$ for $m\in\Z$ we obtain
$$
\lim_{\rho\rightarrow 1-} S_{\Alpha_j, r, \rho}|_{(-\dfrac{n}{2}+l,-\dfrac{n}{2}+l+1)}=
 \sum_{k=0}^{n-1} (l - k)^r \Alpha_j^{l-k} f_k.
$$
The limit on the left is taken in the space of distributions. However, since the limit also exists in the space of continuous functions with uniform convergence on compact sets, and the right hand side is continuous, we also obtain the corresponding identity for the analytic continuation:
$$
S_{\Alpha_j, r}(\Exp{2\pi\I\phi}) = \sum_{k=0}^{n-1} (2\pi\I(l - k))^r \Alpha_j^{l-k} f_k\quad \text{for $\phi\in(-\dfrac{n}2+l,-\dfrac{n}2+l+1)$.}
$$

\pph \coro \label{cor:continuation_a} The analytic continuation of the basis $(S_{\Alpha_j, r}(z))$ with $\dfrac{1}{2\pi}\Arg z \in (-\dfrac{n}2+l,-\dfrac{n}2+l+1)$ is related to the basis $(f_k(\phi))$ by the transposed generalized Vandermonde transformation $\VM_{\Alpha, m_\Alpha, l}^\transpose$: for $\phi\in(-\dfrac{n}2+l,-\dfrac{n}2+l+1)$, $z = \Exp{2\pi\I\phi}$ one has
$$
(S_{\Alpha_j, r}(z)) = (f_k(\phi))\, \VM_{\Alpha, m_\Alpha, l}^\transpose.
$$
\qed
\medskip

The corresponding statement for $S_{\Beta_j, r}(z)$ is completely analogous:

\pph\coro \label{cor:continuation_b} The analytic continuation of the basis $(S_{\Beta_j, r}(z))$ with $\dfrac{1}{2\pi}\Arg z \in (-\dfrac{n}2+l,-\dfrac{n}2+l+1)$ is related to the basis $(f_k(\phi))$ by the transposed generalized Vandermonde transformation $\VM_{\Beta, m_\Beta, l}^\transpose$: for $\phi\in(-\dfrac{n}2+l,-\dfrac{n}2+l+1)$, $z = \Exp{2\pi\I\phi}$ one has
$$
(S_{\Beta_j, r}(z)) = (f_k(\phi))\, \VM_{\Beta, m_\Beta, l}^\transpose.
$$
\qed

\bigskip
In what follows we omit the subscripts $m_\Alpha$, $m_\Beta$, $l$ for the typographic reason.

\pph \label{theorem_on five--tuples} \bf Theorem on hypergeometric monodromy. \rm Choose $\Arg z$ in such a manner that
$$
\dfrac{1}{2\pi} \Arg z \in (-\dfrac{n}{2} + l, -\dfrac{n}{2} + l + 1).
$$
Then the monodromy of the equation $\HH_{\lambda}(\alpha_i, \beta_i)$ with $\lambda = (-1)^n$ is given by the matrices
\benum
\item
$$
M_0^\Alpha = D_\Alpha^\transpose,
$$
\item
$$
M_\infty^\Alpha = \left(\VM_\Alpha {\VM_\Beta}^{-1} {D_\Beta}^{-1} \VM_\Beta {\VM_\Alpha}^{-1}\right)^\transpose \quad \text{in the basis}\; S_{\Alpha_j, r} ,
$$
\item
$$
M_0^\Beta = \left(\VM_\Beta {\VM_\Alpha}^{-1} D_\Alpha \VM_\Alpha {\VM_\Beta}^{-1}\right)^\transpose,
$$
\item
$$
M_\infty^\Beta = \left(D_\Beta^{-1}\right)^\transpose \quad \text{in the basis}\; S_{\Beta_j, r}.
$$
\eenum

\proof The statements ({\rm i}) and ({\rm iv}) are proved in Propositions \ref{prop:monodromy_a} and \ref{prop:monodromy_b}. It follows from Corollaries \ref{cor:continuation_a} and \ref{cor:continuation_b} that
$$
M_0^f = \left({\VM_\Alpha}^{-1} D_\Alpha \VM_\Alpha\right)^\transpose, \quad
M_\infty^f = \left({\VM_\Beta}^{-1} {D_\Beta}^{-1} \VM_\Beta\right)^\transpose \quad \text{in the basis}\; f_k.
$$
The statements ({\rm ii}) and ({\rm iii}) follow.

\pph \rema The monodromy matrices $M_0^f$, $M_\infty^f$ are similar to the generators of the hypergeometric group considered by Levelt \cite{Levelt61} and Beukers and Heckman \cite{BH89}.

\bigskip
\bigskip
\bigskip

\bf Acknowledgements. \rm The first named author thanks Yves Andre and Francesco Baldassarri for the
discussions of the subject. We thank Don Zagier for the reference
to generalized Vandermonde matrices, and Wadim Zudilin for his remarks on the paper.

\bigskip
\bigskip
\bigskip
\bigskip
\bigskip

\nocite{Dwork83}
\nocite{BC04}
\nocite{OTY88}
\nocite{Levelt61}
\nocite{KKP08}
\nocite{BH89}
\nocite{Iritani07}

%
%
%

\end{document}